\documentclass{amsart}
\usepackage{amssymb}
\usepackage{latexsym}
\usepackage{amsmath}
\usepackage{euscript}
\usepackage{graphics}
\usepackage[all]{xy}

\usepackage[leqno]{amsmath}

\usepackage{bbm}

\usepackage{dsfont}

\usepackage{graphicx}

\newcommand{\be}{\begin{equation}}
\newcommand{\ee}{\end{equation}}
\newcommand{\ba}{\begin{eqnarray}}
\newcommand{\ea}{\end{eqnarray}}
\newcommand{\baa}{\begin{eqnarray*}}
\newcommand{\eaa}{\end{eqnarray*}}
\newcommand{\bb}{\begin{thebibliography}}
\newcommand{\eb}{\end{thebibliography}}

\newcommand{\bi}[1]{\bibitem{#1}}
\newcommand{\lab}[1]{\label{#1}}
\newcommand{\re}[1]{(\ref{#1})}

\newcounter{my}
\newcommand{\he}%
   {\stepcounter{equation}\setcounter{my}%
   {\value{equation}}\setcounter{equation}0%
   }%
\newcommand{\she}%
   {\setcounter{equation}{\value{my}}%
    }%

\renewcommand\t{\tilde}

\newcommand\ve{\varepsilon}

\theoremstyle{definition}

\numberwithin{equation}{section}

\begin{document}

\title[Dual -1 Hahn polynomials]
{Dual -1 Hahn polynomials: "classical" polynomials beyond the
Leonard duality}

\author{Satoshi Tsujimoto}
\author{Luc Vinet}
\author{Alexei Zhedanov}

\address{Department of Applied Mathematics and Physics, Graduate School of Informatics,
Kyoto University, Sakyo-ku, Kyoto 606--8501, Japan}

\address{Centre de recherches math\'ematiques,
Universit\'e de Montr\'eal, P.O. Box 6128, Centre-ville Station,
Montr\'eal (Qu\'ebec), H3C 3J7,Canada}

\address{Donetsk Institute for Physics and Technology,\\
83114 Donetsk, Ukraine \\}

%\email{zhedanov@yahoo.com}

%\date{\today}

\begin{abstract}
We introduce the -1 dual Hahn polynomials through  an appropriate
$q \to -1$ limit of the dual q-Hahn polynomials. These polynomials
are orthogonal on a finite set of discrete points on the real
axis, but in contrast to the classical orthogonal polynomials of
the Askey scheme, the -1 dual Hahn polynomials do not exhibit the
Leonard duality property. Instead, these polynomials satisfy a
4-th order difference eigenvalue equation and thus possess a
bispectrality property. The corresponding generalized Leonard pair
consists of two matrices $A,B$ each of size $N+1 \times N+1$. In
the eigenbasis where the matrix $A$ is diagonal, the matrix $B$ is
3-diagonal; but in the eigenbasis where the matrix $B$ is
diagonal, the matrix $A$ is 5-diagonal.
\end{abstract}


\maketitle

\section{Introduction}
\setcounter{equation}{0} Recently new explicit families of
"classical" orthogonal polynomials $P_n(x)$ were introduced
\cite{VZ_little}, \cite{VZ_big}, \cite{VZ_Bochner}. These
polynomials satisfy an eigenvalue equation of the form \be L
P_n(x) = \lambda_n P_n(x). \lab{LPP0} \ee The operator $L$ is of
first order in the derivative operator $\partial_x$ and contains
moreover the reflection operator $R$ defined by $R f(x)=f(-x)$; it
can be identified as a first order operator of Dunkl type written
as \be L = F(x)(I-R) + G(x)\partial_x R  \lab{L_Dunkl} \ee with
some real rational functions $F(x), G(x)$. The corresponding
polynomial eigensolutions $P_n(x)$ can be obtained from the big
and little q-Jacobi polynomials by an appropriate limit $q \to
-1$.

In \cite{BI_TVZ} we generalized this approach to the case of Dunkl
shift operators. In this case the operator $L$ contains the shift
operator $T^+ f(x) =f(x+1)$ together with the reflection operator
$R$: \be L = F(x)(I-R) + G(x)(T^+ R-I) \lab{L_FG} \ee ($I$ stands
for the identity operator). The rational functions can be
recovered from the condition that the operator $L$ stabilizes the
spaces of polynomials, i.e. it sends any polynomial of degree $n$
into a polynomial of the same degree. It can then be demonstrated
\cite{BI_TVZ} that the polynomial eigensolutions $P_n(x)$ of the
eigenvalue equation \re{LPP0} with the operator \re{L_FG} satisfy
the 3-term recurrence relation and hence are orthogonal
polynomials. In fact, the polynomials $P_n(x)$ in this case
coincide with the Bannai-Ito (BI) polynomials first constructed in
\cite{BI} (see also \cite{Ter}).

The BI polynomials thus possess  the bispectrality property: they
satisfy simultaneously the 3-term recurrence relation (common to
all orthogonal polynomials) \be P_{n+1}(x) + b_n P_n(x) + u_n
P_{n-1}(x) = xP_n(x) \lab{rec_P_gen} \ee and the eigenvalue
equation \re{LPP0}. Moreover, in the case when the support of the
orthogonality measure consists of a finite number of points $x_s,
\: s=0,1,2,\dots,N$, the BI polynomials satisfy the Leonard
duality property \cite{Leonard}, \cite{BI}, \cite{Ter}. This means
that there is a finite difference equation of the form \be U_s
\left(P_n(x_{s+1}) - P_n(x_s) \right) + V_s \left(P_n(x_{s-1}) -
P_n(x_s) \right) = \lambda_n P_n(x_s) \lab{Leon_BI} \ee with some
real coefficients $V_s, U_s$. The difference equation \re{Leon_BI}
is of the second order and can be considered as a dual relation
with respect to the recurrence relation. In fact, the difference
equation equation \re{Leon_BI} is a simple consequence of the
eigenvalue equation \re{LPP0} with the Dunkl shift operator
\re{L_FG} \cite{BI_TVZ}.

We showed in \cite{BI_TVZ} that the BI polynomials can be obtained
by an appropriate $q \to -1$ limit from the Askey-Wilson
polynomials. Correspondingly, the Dunkl shift operator $L$ appears
in the same limit from the Askey-Wilson difference operator
\cite{BI_TVZ}. It should be stressed that there are several
possibilities in taking the limit $q \to -1$ of the Askey-Wilson
polynomials. Not all of them lead to the BI polynomials. There is
another family of orthogonal polynomials - called the
complementary BI polynomials in \cite{BI_TVZ} - which can be
obtained from the Askey-Wilson polynomials in the same limit. In
contrast to the BI polynomials, the complementary BI (CBI)
polynomials do not satisfy an eigenvalue equation of the form
\re{LPP0}. As a consequence, the CBI polynomials do not possess
the Leonard duality property \re{Leon_BI}.

The main purpose of the present paper is to study the orthogonal
polynomials which appear in the $q \to -1$ limit of the dual
q-Hahn polynomials. We call these polynomials the dual -1 Hahn
polynomials. We derive explicitly basic properties and relations
for them including a 3-term recurrence relation. The main result
is the existence of a dual eigenvalue equation \re{LPP0} for the
dual -1 Hahn polynomials. In contrast to the BI polynomials, the
operator $L$ is now 2-nd order with respect to the shift
operators. As a consequence, the dual -1 Hahn polynomials obey a
5-term difference relation on an appropriate grid $x_s$ (instead
of the 3-term relation \re{Leon_BI} for the BI polynomials).

\section{Dual q-Hahn polynomials}
\setcounter{equation}{0} The dual $q$-Hahn polynomials \cite{KLS}
\be R_n(x;a,b,N) = {_3}\Phi_2 \left( {q^{-n}, q^{-s}, abq^{s+1}
\atop aq, q^{-N} } \left | \right . q;q  \right) \lab{dual_q_Hahn}
\ee depend on 3 parameters $a,b,N$, where $N$ is a positive
integer. The argument $x$ can be parametrized in terms of $s$ as
follows \be x = q^{-s} + abq^{s+1}. \lab{x_s} \ee The polynomials
$R_n(x;a,b,N)$ satisfy the 3-term recurrence relation \be A_n
R_{n+1}(x;a,b,N) - (A_n+C_n-1-abq) R_{n}(x;a,b,N) + C_n
R_{n-1}(x;a,b,N) = xR_{n}(x;a,b,N), \lab{3_term_Hahn} \ee where
\be A_n = (1-q^{n-N})(1-aq^{n+1}), \quad C_n =
aq(1-q^n)(b-q^{n-N-1}). \lab{AC-q_Hahn} \ee The polynomials
$R_n(x;a,b,N)$ are not monic; their monic version $P_n(x;a,b,N)=
\kappa_n R_n(x;a,b,N) =x^n + O(x^{n-1})$ (with $\kappa_n$ an
appropriate factor) will obey the recurrence relation \be
P_{n+1}(x;a,b,N) + b_n P_n(x;a,b,N) + u_n P_{n-1}(x;a,b,N) =
xP_{n}(x;a,b,N), \lab{rec_monic_q} \ee where
$$
b_n = 1+abq-A_n-C_n, \quad u_n = A_{n-1} C_n .
$$
The dual q-Hahn polynomials also verify a q-difference equation of
second order \cite{KLS} \be B(s) R_n(x_{s+1}) +D(s)R_n(x_{s-1})
-(B(s)+D(s))R_n(x_s) =(q^{-n}-1)R_n(x_s), \quad s=0,1,2,\dots, N,
\lab{dfrc_q} \ee where $x_s$ is given by \re{x_s} and \ba &&B(s) =
\frac{(1-q^{s-N})(1-aq^{s+1})(1-abq^{s+1})}{(1-abq^{2s+1})(1-abq^{2s+2})}
\nonumber \\ &&D(s) =
-\frac{aq^{s-N}(1-q^{s})(1-abq^{s+1+N})(1-bq^{s})}{(1-abq^{2s+1})(1-abq^{2s})}.
\lab{BD_q} \ea As per \cite{NSU}, the equation \re{dfrc_q} can be
interpreted as a q-difference equation on the "q-quadratic grid"
$x_s$.

The dual q-Hahn polynomials satisfy the orthogonality relation \be
\sum_{s=0}^N w_s R_n(x_s)R_m(x_s) = \kappa_0 u_1 u_1 \dots u_n
\delta_{nm}, \lab{ort_DH} \ee where the discrete weights are \be
w_s =\frac{(aq,abq,q^{-N};q)_s}{(q,abq^{N+2},bq;q)_s}
\frac{1-abq^{2s+1}}{(1-abq)(-aq)^s}\: q^{Ns-s(s-1)/2} \lab{w_qH}
\ee and the normalization constant is
$$
\kappa_0=\frac{(abq^2;q)_N}{(bq;q)_N} (aq)^{-N}.
$$
We used the standard notation for the q-shifted factorials
\cite{KLS}.

When $q \to 1$ and $a=q^{\alpha}, \: b=q^{\beta}$ the dual q-Hahn
polynomials become the ordinary dual Hahn polynomials \cite{KLS}
\be W_n(x_s;a,b,N) = {_3}F_2 \left({ -n, -s, s+1+\alpha+\beta
\atop \alpha+1, -N } \left | \right . 1 \right) \lab{dual_Hahn}
\ee where
$$
x_s=s(s+\alpha+\beta+1) .
$$
These polynomials satisfy the three-term recurrence relation \be
A_n W_{n+1}(x_s) - (A_n+C_n)W_n(x_s) + C_n W_{n-1}(x_s)  = x_s
W_n(x_s), \lab{rec_W} \ee where
$$
A_n = (n-N)(n+\alpha+1), \; C_n = n(n-\beta-N-1).
$$
The corresponding monic dual Hahn polynomials $\hat W_n(x)$ obey
\be \hat W_{n+1}(x_s)+ b_n \hat W_{n}(x_s) + u_n \hat W_{n-1}(x_s)
= x_s \hat W_{n}(x_s), \lab{monic_W_rec} \ee where \be u_n =
A_{n-1}C_n = n(n-\beta-N-1)(n-N-1)(n+\alpha), \quad b_n =
-A_n-C_n. \lab{ub_W} \ee

\section{A limit $q \to -1$ of the recurrence relations and the orthogonality property}
\setcounter{equation}{0} We wish to consider a limit $q \to -1$.
We will assume that $a \to \pm 1$ and $b \to \pm 1$ when $q \to
-1$. We want to obtain a nondegenerate limit of the coefficients
$A_n(1+q)^{-1}, C_n(1+q)^{-1}$ for $q \to -1$. This means that
both these limit coefficients should exist and be nonzero for all
admissible values $n=1,2, \dots, N$. It is hence easily seen that
necessarily, we must have $ab \to 1$. Two situations have to be
considered separately:

(i) when $N=2,4,6,\dots$ is even, the nontrivial $q \to -1$ limit
then  exists iff $a \to 1, \; b \to 1$. It is then natural to take
the parametrization \be q=-e^{\ve}, \; a=e^{-\alpha \ve}, \; b =
e^{-\beta \ve}, \quad \ve \to 0 \lab{ab_even} \ee with real
parameters $\alpha,\beta$. Dividing the recurrence relation
\re{3_term_Hahn} by $q+1$ and taking the limit $\ve \to 0$ we
obtain the recurrence relation \be A_n^{(-1)}
R_{n+1}^{(-1)}(y_s;\alpha,\beta,N) + C_n^{(-1)}
R_{n-1}^{(-1)}(y_s;\alpha,\beta,N)- \lab{rec_N_even} \ee
$$
(A_n^{(-1)}+ C_n^{(-1)}) R_{n}^{(-1)}(y_s;\alpha,\beta,N))= y_s
R_n^{(-1)}(y_s;\alpha,\beta,N),  $$ where the grid $y_s$ has the
following expression \be y_s=\left\{ { -\alpha-\beta+2s+1 \quad
\mbox{if} \quad s \quad \mbox{even}, \atop \alpha+\beta-2s-1 \quad
\mbox{if} \quad s \quad \mbox{odd}} \right . \lab{even_grid_y} \ee
or, equivalently, \be y_s=(-1)^s \: (1-\alpha-\beta+2s), \quad
s=0,1,\dots, N. \lab{y_s_N_even} \ee The recurrence coefficients
are \be A_n^{(-1)}=\left\{ { 2(n-N) \quad \mbox{if} \quad n \quad
\mbox{even}   \atop 2(n+1-\alpha) \quad \mbox{if} \quad n \quad
\mbox{odd}}  \right . , \quad C_n^{(-1)}=\left\{ { -2n \quad
\mbox{if} \quad n \quad \mbox{even} \atop 2(N+1-\beta -n) \quad
\mbox{if} \quad n \quad \mbox{odd}} \right . . \lab{AC_even} \ee
The corresponding monic $-1$ dual Hahn polynomials satisfy
relation \re{rec_monic_q}, where \be
u_n^{(-1)}=A_{n-1}^{(-1)}C_n^{(-1)}=\left\{ { 4n(\alpha-n) \quad
\mbox{if} \quad n \quad \mbox{even}   \atop 4(N-n+1)(n+\beta-N-1)
\quad \mbox{if} \quad n \quad \mbox{odd}}  \right . \lab{u-1_even}
\ee and \be b_n^{(-1)}=1-\alpha-\beta-A_{n}^{(-1)}
-C_n^{(-1)}=\left\{ {2N+1-\alpha-\beta \quad \mbox{if} \quad n
\quad \mbox{even}   \atop -2N-3 +\alpha+\beta \quad \mbox{if}
\quad n \quad \mbox{odd}}  \right . \lab{b-1_even} \ee It is
convenient to introduce the "$\mu$-numbers" \be [n]_{\mu} =
n+\mu(1-(-1)^n), \lab{mu_num} \ee which appear naturally in
problems connected with the Dunkl operators \cite{Rosen}. One can
then present the recurrence coefficients in the compact form \be
u_n^{(-1)} = 4 [n]_{\xi} [N-n+1]_{\eta}, \quad b_n^{(-1)}=
2([n]_{\xi} + [N-n]_{\eta}) + 1-\alpha-\beta, \lab{rec_mu_form_e}
\ee where
$$
\xi=\frac{\beta-N-1}{2}, \; \eta=\frac{\alpha-N-1}{2}.
$$

It is seen that $u_0=u_{N+1}=0$ as required for finite orthogonal
polynomials. The positivity condition $u_n>0, \; n=1,2,\dots, N$
is equivalent to the conditions \be \alpha>N, \; \beta>N.
\lab{pos_even_cond} \ee

These polynomials are orthogonal on the finite set of points $y_s$
\be \sum_{s=0}^N w_s R_{n}^{(-1)}(y_s) R_{m}^{(-1)}(y_s) =
\kappa_0 u_1^{(-1)} u_2^{(-1)} \dots u_n^{(-1)} \: \delta_{nm},
\lab{ort_N_even} \ee where the discrete weights are defined as \be
w_{2s} = (-1)^s \frac{(-N/2)_s}{s!} \:
\frac{(1-\alpha/2)_s(1-\alpha/2-\beta/2)_s}{(1-\beta/2)_s
(N/2+1-\alpha/2-\beta/2)_s}, \quad s=0,1,2,\dots, {N\over 2}
\lab{w_N_e_s_e} \ee and \be w_{2s+1} = (-1)^s
\frac{(-N/2)_{s+1}}{s!} \:
\frac{(1-\alpha/2)_s(1-\alpha/2-\beta/2)_s}{(1-\beta/2)_s
(N/2+1-\alpha/2-\beta/2)_{s+1}}, \quad s=0,1,\dots, {N\over 2}-1.
\lab{w_N_e_s_o} \ee The normalization coefficient is \be \kappa_0
= \frac{\left( 1-\frac{\alpha+\beta}{2} \right)_{N/2 }}{\left(
1-\frac{\beta}{2} \right)_{N/2 }}. \lab{kap_even} \ee

Assume that $\alpha=N+\epsilon_1, \; \beta=N+\epsilon_2$, where
$\epsilon_{1,2}$ are arbitrary positive parameters. This
parametrization corresponds to the positive condition for the dual
-1 Hahn polynomials. Then it is easily verified that all the
weights are positive $w_s>0, \: s=0,1,\dots,N$.

Moreover, the spectral points $y_s$ are divided into two
non-overlapping discrete sets of the real line: $$\{1-\delta,
-3-\delta, -7-\delta, \dots, -2N+1-\delta \}$$ and
$$
\{1+\delta, 5+\delta, 9+\delta, \dots, 2N-3+\delta \},
$$
where $\delta=\epsilon_1+\epsilon_2>0$. The first set corresponds
to $y_s$ with even $s$ and contains $1+N/2$ points; the second set
corresponds to $y_s$ with odd $s$ and contains $N/2$ points.

(ii) when $N=1,3,5,\dots$ is odd, a nontrivial $q \to -1$ limit
also exists iff $a \to -1, \; b \to -1$. We take the
parametrization \be q=-e^{\ve}, \; a=-e^{\alpha \ve}, \; b =
-e^{\beta \ve}, \quad \ve \to 0 \lab{ab_odd} \ee with real
parameters $\alpha,\beta$. Dividing again the recurrence relation
\re{3_term_Hahn} by $q+1$ and taking the limit $\ve \to 0$, we
obtain the recurrence relation \re{rec_N_even} where the grid
$y_s$ is defined as \be y_s=\left\{ { \alpha+\beta+2s+1 \quad
\mbox{if} \quad s \quad \mbox{even} \atop -\alpha-\beta-2s-1 \quad
\mbox{if} \quad s \quad \mbox{odd}} \right . \lab{odd_grid_y} \ee
and the recurrence coefficients given by \be A_n^{(-1)}n=\left\{ {
2(\alpha+n+1) \quad \mbox{if} \quad n \quad \mbox{even}   \atop
2(n-N) \quad \mbox{if} \quad n \quad \mbox{odd}}  \right . , \quad
C_n^{(-1)}=\left\{ { -2n \quad \mbox{if} \quad n \quad \mbox{even}
\atop 2(\beta+N -n+1) \quad \mbox{if} \quad n \quad \mbox{odd}}
\right .  .   \lab{AC_odd}  \ee The corresponding monic $-1$ dual
Hahn polynomials satisfy the standard relation \re{rec_monic_q},
where \be u_n^{(-1)}=A_{n-1}^{(-1)} C_n^{(-1)}=\left\{ { 4n(N+1-n)
\quad \mbox{if} \quad n \quad \mbox{even}   \atop
4(\alpha+n)(\beta+N+1-n) \quad \mbox{if} \quad n \quad \mbox{odd}}
\right . \lab{u-1_odd} \ee and \be
b_n^{(-1)}=1+\alpha+\beta-A_{n}^{(-1)} -C_n^{(-1)}=\left\{
{-1-\alpha+\beta \quad \mbox{if} \quad n \quad \mbox{even}   \atop
-1 +\alpha-\beta \quad \mbox{if} \quad n \quad \mbox{odd}}  \right
. . \lab{b-1_odd} \ee Again, as in the case of even $N$, it is
possible to present the recurrence coefficients in the compact
form \be u_n^{(-1)} = 4 [n]_{\xi} [N-n+1]_{\eta}, \quad
b_n^{(-1)}= 2([n]_{\xi} + [N-n]_{\eta}) -2N-1-\alpha-\beta,
\lab{rec_mu_form_o} \ee with $\xi=\alpha/2, \: \eta=\beta/2$. It
is seen that in both cases: $N$ even and $N$ odd, the recurrence
coefficients of the dual -1 Hahn polynomials are presented in the
unified form \re{rec_mu_form_e} or \re{rec_mu_form_o} with the
difference only residing with the parameters $\xi,\eta$.

It is seen that $u_0=u_{N+1}=0$ as required for finite orthogonal
polynomials. The positivity condition $u_n>0, \; n=1,2,\dots, N$
is equivalent either to condition \be \alpha>-1, \; \beta>-1
\lab{pos_odd_cond} \ee or to condition $\alpha<-N, \: \beta<-N$.
In what follows we shall use only condition \re{pos_odd_cond}.

The polynomials $R_n^{(-1)}(x)$ are orthogonal on the finite set
of points $y_s$ \be \sum_{s=0}^N w_s R_{n}^{(-1)}(y_s)
R_{m}^{(-1)}(y_s) = \kappa_0 u_1^{(-1)} u_2^{(-1)} \dots
u_n^{(-1)} \: \delta_{nm}, \lab{ort_N_odd} \ee where the discrete
weights are defined as \be w_{2s} = (-1)^s \frac{(-(N-1)/2)_s}{s!}
\: \frac{(1/2+\alpha/2)_s(1+\alpha/2+\beta/2)_s}{(1/2+\beta/2)_s
(N/2+3/2+\alpha/2+\beta/2)_s}, \quad s=0,1,2,\dots, {{N-1}\over 2}
\lab{w_O_e_s_e} \ee and \be w_{2s+1} = (-1)^s
\frac{(-(N-1)/2)_s}{s!} \:
\frac{(1/2+\alpha/2)_{s+1}(1+\alpha/2+\beta/2)_s}{(1/2+\beta/2)_{s+1}
(N/2+3/2+\alpha/2+\beta/2)_s}, \quad s=0,1,2,\dots, {{N-1}\over 2}
\lab{w_O_e_s_o} \ee The normalization coefficient is \be \kappa_0
= \frac{\left( 1+\frac{\alpha+\beta}{2} \right)_{(N+1)/2 }}{\left(
\frac{\beta+1}{2} \right)_{(N+1)/2 }}. \lab{kap_odd} \ee

Assume that $\alpha=-1+\epsilon_1, \; \beta=-1+\epsilon_2$, where
$\epsilon_{1,2}$ are arbitrary positive parameters. This
parametrization corresponds to the positive condition for the dual
-1 Hahn polynomials for $N$ odd. Then it is easily verified that
the weights are positive $w_s>0, \: s=0,1,\dots,N$.

Moreover, the spectral points $y_s$ are divided into two
non-overlapped discrete sets of the real line: $$\{-1-\delta,
-5-\delta, -9-\delta, \dots, -2N+1-\delta \}$$ and
$$
\{-1+\delta, 3+\delta, 7+\delta, \dots, 2N-3+\delta \},
$$
where $\delta=\epsilon_1+\epsilon_2>0$. Both sets contain
$(N-1)/2$ points.

\section{Explicit expression in terms of the ordinary dual Hahn polynomials}
\setcounter{equation}{0} In this section, we derive an explicit
expression for the -1 dual Hahn polynomials in terms of the
ordinary dual Hahn polynomials.

Consider first the case of even $N=2,4,6,\dots$. Introduce the
"shifted" monic -1 dual Hahn polynomials $\tilde R_n(x)=
R^{(-1)}_n(x-1)$ From formulas \re{u-1_even}, \re{b-1_even}, we
conclude that these polynomials satisfy the recurrence relation
\be \t R_{n+1}(x) + (-1)^n \tau \t R_{n}(x) + u_n \t R_{n-1}(x) =
x \t R_{n}(x), \lab{rec_shift_even} \ee where
$$
\tau=2N+2-\alpha-\beta
$$
and $u_n$ are given by \re{u-1_even}.

A recurrence relation of the type \re{rec_shift_even} leads to
orthogonal polynomials $\t R_n(x)$ which are very close to
symmetric orthogonal polynomials. Using methods developed in
\cite{Chi_Bol} and \cite{VZ_big}, we can introduce a pair of monic
orthogonal polynomials $P_n(x)$ and $Q_n(x)$ by the formulas : \be
\t R_{2n}(x) =P_n(x^2), \quad \t R_{2n+1}(x)=(x-\tau)Q_n(x^2) .
\lab{RPQ} \ee It can easily be shown that the polynomials $P_n(x)$
and $Q_n(x)$ satisfy the following recurrence relations (it is
assumed that $u_0=0$) \be P_{n+1}(x) +(u_{2n} + u_{2n+1}+
\tau^2)P_n(x)+ u_{2n} u_{2n-1} P_{n-1}(x)=xP_n(x) \lab{rec_P} \ee
and \be Q_{n+1}(x) +(u_{2n+2} + u_{2n+1}+ \tau^2)Q_n(x)+ u_{2n}
u_{2n+1} Q_{n-1}(x)=xQ_n(x), \lab{rec_Q} \ee and moreover that the
polynomials are connected by the Christoffel transform \be Q_n(x)
= \frac{P_{n+1}(x) + u_{2n+1}P_n(x)}{x-\tau^2}. \lab{PQ_Chr} \ee
It is also easily seen that both $P_n(x)$ and $Q_n(x)$ are
ordinary dual Hahn polynomials. We hence have the following
explicit expression: \be R_{2n}^{(-1)}(x-1) = \gamma^{(0)}_n \:
{_3}F_{2} \left( {-n, \eta+\frac{x}{4}, \eta-\frac{x}{4}  \atop
-\frac{N}{2}, 1-\frac{\alpha}{2}}   \left | \right . 1 \right),
\quad n=0,1,2,\dots \lab{R_even_exp} \ee and \be
R_{2n+1}^{(-1)}(x-1) = \gamma^{(1)}_n \: (x-\tau) \: {_3}F_{2}
\left( {-n, \eta+\frac{x}{4}, \eta-\frac{x}{4}  \atop
1-\frac{N}{2}, 1-\frac{\alpha}{2}}   \left | \right . 1 \right),
\quad n=0,1,2,\dots , \lab{R_odd_exp} \ee where
$\eta=1/2-(\alpha+\beta)/4$ and the normalization coefficients are
$$
\gamma^{(0)}_n= 16^n (-N/2)_n(1-\alpha/2)_n, \quad \gamma^{(1)}_n=
16^n (1-N/2)_n(1-\alpha/2)_n .
$$
Quite similarly, for the odd $N=1,3,5,\dots$ we find \be
R_{2n}^{(-1)}(x-1) = \gamma^{(0)}_n \: {_3}F_{2} \left( {-n,
\eta+\frac{x}{4}, \eta-\frac{x}{4}  \atop -\frac{N-1}{2},
\frac{\alpha+1}{2}}   \left | \right . 1 \right), \quad
n=0,1,2,\dots \lab{RO_even_exp} \ee and \be R_{2n}^{(-1)}(x-1) =
\gamma^{(1)}_n (x+\alpha-\beta) \: {_3}F_{2} \left( {-n,
\eta+\frac{x}{4}, \eta-\frac{x}{4}  \atop -\frac{N-1}{2},
\frac{\alpha+3}{2}}   \left | \right . 1 \right), \quad
n=0,1,2,\dots, \lab{RO_odd_exp} \ee where
$$
\eta=\frac{\alpha+\beta+2}{4}, \; \gamma_n^{(0)}=16^n
\left(\frac{1-N}{2}\right)_n \left(\frac{\alpha+1}{2}\right)_n, \;
\gamma_n^{(1)}=16^n \left(\frac{1-N}{2}\right)_n
\left(\frac{\alpha+3}{2}\right)_n .
$$
Some of these polynomials have appeared in \cite{Stoilova},
\cite{JSJ} in the context of quantum spin chains.

\section{Difference equation}
\setcounter{equation}{0} Consider the following operator $L$
defined on the space of functions $f(s)$ that depend on a discrete
variable $s$: \be L f(s) = B(s) f(s+1) + D(s) f(s-1) - (B(s)
+D(s))f(s), \quad s=0,1,2,\dots ,  \lab{L_q_def} \ee where $B(s),
D(s)$ are given in \re{BD_q}. Manifestly, the difference equation
\re{dfrc_q} means that the dual $q$-Hahn polynomials are
eigenfunctions of the operator $L$ \be L R_n(x(s)) = \lambda_n
R_n(x(s)) \lab{LR_eig} \ee with the eigenvalues \be \lambda_n =
q^{-n}-1 . \lab{lambda_q} \ee When $q \to 1$ we obtain the
difference eigenvalue equation for the ordinary dual Hahn
polynomials \be L_1 W_n(x(s)) = -n W_n(x(s)), \lab{LR_eig_1} \ee
where the operator $L_1$ can be obtained from $L$ as $L_1 =
\lim_{q \to 1} L(q-1)^{-1}$. Explicitly \cite{KLS} \be L_1 =
B_1(s) f(s+1) + D_1(s) f(s-1) - (B_1(s) +D_1(s))f(s),
\lab{L_1_Hahn} \ee where
$$
B_1(s) =
\frac{(s+\alpha+\beta+1)(s+\alpha+1)(N-s)}{(2s+\alpha+\beta+1)(2s+\alpha+\beta+2)},
\quad D_1(s) =
\frac{s(s+\beta)(s+\alpha+\beta+N+1)}{(2s+\alpha+\beta+1)(2s+\alpha+\beta)}
.
$$
When we try to perform a similar procedure for the limit $q \to
-1$ we encounter a problem. Indeed, it is easily seen that
$L(1+q)^{-1}$ does not have a nondegenerate limit as $q \to -1$.
We thus cannot obtain an eigenvalue equation in 3-diagonal form
like \re{LR_eig_1} for the dual -1 Hahn polynomials.

Nevertheless we observe that the operator $(L^2 + 2L)(1+q)^{-1}$
does survive in the limit $q \to -1$. We hence have the following
eigenvalue equation for the dual -1 Hahn polynomials \be H
R_n^{(-1)}(y_{s}) = 2n R_n^{(-1)}(y_s), \lab{HR_eig} \ee where the
grid $y_s$ is defined by \re{even_grid_y} or \re{odd_grid_y} and
$$
H = \lim_{q \to -1}(L^2 + 2L)(1+q)^{-1}.
$$
The operator $H$ obviously is 5-diagonal, i.e. \be H f(s)=
U_2(s)(f(s+2)-f(s)) + U_1(s)(f(s+1)-f(s)) + V_2(s)(f(s-2)-f(s)) +
V_1(s)(f(s-1)-f(s)) \lab{H-5-diag} \ee The explicit expressions
for the coefficients $U_i(s), V_i(s)$ depend on the parity of $N$.
For even $N=2,4,6,\dots$ they are: \be U_2(s)=\left\{ { -2\,{\frac
{ \left( \alpha-s-2 \right) \left( \beta+\alpha-s-2
 \right)  \left( N-s \right) }{ \left( \alpha+\beta-2\,s-2 \right)
 \left( \alpha+\beta-2\,s-4 \right) }}
 \quad \mbox{if} \quad s \quad \mbox{even}   \atop -2\,{\frac { \left( \alpha-s-1 \right)  \left( \beta+\alpha-s-1
 \right)  \left( N-s-1 \right) }{ \left( \alpha+\beta-2\,s-2 \right)
 \left( \alpha+\beta-2\,s-4 \right) }}
 \quad \mbox{if} \quad s \quad \mbox{odd}}  \right . , \lab{U2_ev} \ee
\be
U_1(s)=\left\{ { 2\,{\frac { \left( \beta+\alpha \right)  \left( \alpha-\beta \right)
 \left( N-s \right) }{ \left( \alpha+\beta-2\,s \right)  \left( \alpha
+\beta-2\,s-2 \right)  \left( \alpha+\beta-2\,s-4 \right) }}\quad \mbox{if} \quad s \quad \mbox{even}   \atop 4\,{\frac { \left( \alpha-s-1 \right)  \left( \beta+\alpha-s-1
 \right)  \left( 2\,N+2-\alpha-\beta \right) }{ \left( \alpha+\beta-2
\,s \right)  \left( \alpha+\beta-2\,s-2 \right)  \left(
\alpha+\beta-2 \,s-4 \right) }}\quad \mbox{if} \quad s \quad
\mbox{odd}}  \right . , \lab{U1_ev} \ee \be V_2(s)=\left\{
{2\,{\frac {s \left( \beta-s \right)  \left( -\alpha-\beta+N+s
\right) }{ \left( \alpha+\beta-2\,s+2 \right)  \left(
\alpha+\beta-2\,s
 \right) }}\quad \mbox{if} \quad s \quad \mbox{even}   \atop 2\,{\frac { \left( s-1 \right)  \left( \beta-s+1 \right)  \left( -
\alpha-\beta+N+s+1 \right) }{ \left( \alpha+\beta-2\,s+2 \right)
 \left( \alpha+\beta-2\,s \right) }}\quad \mbox{if} \quad s \quad \mbox{odd}}  \right . , \lab{V2_ev} \ee
\be
V_1(s)=\left\{ {4\,{\frac {s \left( \beta-s \right)  \left( 2\,N+2-\alpha-\beta
 \right) }{ \left( \alpha+\beta-2\,s \right)  \left( \alpha+\beta-2\,s
-2 \right)  \left( \alpha+\beta-2\,s+2 \right) }}
\quad \mbox{if} \quad s \quad \mbox{even}   \atop -2\,{\frac { \left( \beta+\alpha \right)  \left( \alpha-\beta \right)
 \left( -\alpha-\beta+N+s+1 \right) }{ \left( \alpha+\beta-2\,s
 \right)  \left( \alpha+\beta-2\,s-2 \right)  \left( \alpha+\beta-2\,s
+2 \right) }}\quad \mbox{if} \quad s \quad \mbox{odd}}  \right . ;
\lab{V1_ev} \ee and for odd $N=1,3,5,\dots$  \be U_2(s)=\left\{ {
-2\,{\frac { \left( \alpha+\beta+s+2 \right) \left( \alpha+s+1
 \right)  \left( N-s-1 \right) }{ \left( \alpha+\beta+2\,s+2 \right)
 \left( \alpha+\beta+2\,s+4 \right) }}\quad \mbox{if} \quad s \quad \mbox{even}   \atop -2\,{\frac { \left( \alpha+s+2 \right)  \left( \alpha+\beta+s+1
 \right)  \left( N-s \right) }{ \left( \alpha+\beta+2\,s+2 \right)
 \left( \alpha+\beta+2\,s+4 \right) }}\quad \mbox{if} \quad s \quad \mbox{odd}}  \right . ,\lab{U2_od} \ee
\be U_1(s)=\left\{ { -2\,{\frac { \left( \alpha+\beta \right)
\left( \alpha+s+1 \right)
 \left( \alpha+\beta+2\,N+2 \right) }{ \left( \alpha+\beta+2\,s
 \right)  \left( \alpha+\beta+2\,s+2 \right)  \left( \alpha+\beta+2\,s
+4 \right) }} \quad \mbox{if} \quad s \quad \mbox{even}   \atop
-4\,{\frac { \left( \alpha-\beta \right)  \left( N-s \right)
\left( \alpha+\beta+s+1 \right) }{ \left( \alpha+\beta+2\,s
\right)  \left( \alpha+\beta+2\,s+2 \right)  \left(
\alpha+\beta+2\,s+4 \right) }} \quad \mbox{if} \quad s \quad
\mbox{odd}}  \right . , \lab{U1_od} \ee \be V_2(s)=\left\{
{-2\,{\frac {s \left( \beta+s-1 \right)  \left( \alpha+\beta+N+s+1
 \right) }{ \left( \alpha+\beta+2\,s-2 \right)  \left( \alpha+\beta+2
\,s \right) }}\quad \mbox{if} \quad s \quad \mbox{even}   \atop
-2\,{\frac { \left( s-1 \right)  \left( \beta+s \right)  \left(
\alpha +\beta+N+s \right) }{ \left( \alpha+\beta+2\,s-2 \right)
\left( \alpha+\beta+2\,s \right) }}\quad \mbox{if} \quad s \quad
\mbox{odd}}  \right . , \lab{V2_od} \ee \be V_1(s)=\left\{
{-4\,{\frac {s \left( \alpha-\beta \right) \left(
\alpha+\beta+N+s+1
 \right) }{ \left( \alpha+\beta+2\,s \right)  \left( \alpha+\beta+2\,s
+2 \right)  \left( \alpha+\beta+2\,s-2 \right) }}\quad \mbox{if} \quad s \quad \mbox{even}   \atop -2\,{\frac { \left( \beta+s \right)  \left( \alpha+\beta \right)
 \left( \alpha+\beta+2\,N+2 \right) }{ \left( \alpha+\beta+2\,s
 \right)  \left( \alpha+\beta+2\,s+2 \right)  \left( \alpha+\beta+2\,s
-2 \right) }}\quad \mbox{if} \quad s \quad \mbox{odd}}  \right . .
\lab{V1_od} \ee

We thus have a difference equation for the dual -1 Hahn
polynomials in the form \ba
&&U_2(s)\left(R_n^{(-1)}(y_{s+2})-R_n^{(-1)}(y_{s}) \right) + U_1(s)\left(R_n^{(-1)}(y_{s+1})-R_n^{(-1)}(y_{s})\right ) +  \nonumber \\
&&V_2(s)\left(R_n^{(-1)}(y_{s-2})-R_n^{(-1)}(y_{s}) \right) +
V_1(s) \left(R_n^{(-1)}(y_{s-1})-R_n^{(-1)}(y_{s}) \right) = 2n
R_n^{(-1)}(y_s), \nonumber \ea where the coefficients $U_i(s),
V_i(s)$ are provided in the above formulas.

\section{Another form of the difference equation}
\setcounter{equation}{0} The difference equation for the dual -1
Hahn polynomials can be presented in a more compact form if one
notices that the grid $y_s$ satisfies the relations \be y_{s\pm
1}=\left\{ {-y_s \mp 2 \quad \mbox{if} \quad s \quad \mbox{even}
\atop   -y_s \pm 2   \quad \mbox{if} \quad s \quad \mbox{odd}}
\right . . \lab{grid_shift} \ee This property implies that the
difference equation can be written as \ba &&E_1(x)
\left(R^{(-1)}_{n}(x+4) -  R^{(-1)}_{n}(x) \right) + E_2(x)
\left(R^{(-1)}_{n}(x-4) -  R^{(-1)}_{n}(x) \right) + \nonumber \\
\nonumber \\ &&G_1(x) \left(R^{(-1)}_{n}(-x-2) - R^{(-1)}_{n}(x)
\right) + G_2(x) \left(R^{(-1)}_{n}(-x+2) - R^{(-1)}_{n}(x)
\right)= 2n R_n^{(-1)}(x) \nonumber \ea or, in operator form as
\be H R_n^{(-1)}(x) = 2n R_n^{(-1)}(x). \lab{HR_op} \ee The
operator $H$ in \re{HR_op} reads \be H=E_1(x) T^4 + E_2(x) T^{-4}
+ G_1(x) T^2R + G_2(x) T^{-2}R -(E_1(x)+E_2(x) + G_1(x) +
G_2(x))I, \lab{op_H} \ee where the operators $T$ and $T^{-1}$ are
the standard shift operators: $T^jf(x) = f(x+j), \: j=0, \pm 1,
\pm 2, \dots$, and $R$ is the reflection operator $Rf(x) = f(-x)$.
$I$ denotes the identity operator.

The functions $E_i(x), G_i(x), \; i=1,2$ are simple rational
functions in $x$.

For even $N=2,4,\dots$ we have
$$
E_1(x) = {\frac { \left(x+3 -\alpha+\beta \right)  \left( x+3-\alpha-\beta
 \right)  \left( x-1-2\,N+\alpha+\beta \right) }{4 \left( x+1
 \right)  \left( 3+x \right) }},$$
$$
E_2(x) = -{\frac { \left( \alpha-1+\beta+x \right)  \left( \alpha-1-\beta+
x \right)  \left( x-1+2\,N-\alpha-\beta \right) }{ 4\left( x-1
 \right)  \left( x-3 \right) }},
$$
$$
G_1(x)={\frac { \left( \alpha^2-\beta^2 \right)
 \left( x+\alpha+\beta -2\,N-1 \right) }{ \left(x^2-1 \right)  \left(
x+3 \right)  }},
$$
$$
G_2(x)={\frac { \left( 2N+2-\alpha-\beta \right) \left( (x+\alpha-1)^2-\beta^2 \right)   }{ \left( x^2-1 \right)
 \left( x-3 \right)   }}
$$
For odd $N=1,3,5,\dots$
$$
E_1(x) = {\frac { \left( x+\alpha+\beta+3 \right)  \left( x+\alpha-\beta+
1 \right)  \left( x-\alpha-\beta-2\,N +1\right) }{4 \left( x+1
 \right)  \left( x+3 \right) }}
$$
$$
E_2(x) = -{\frac { \left( x-\alpha-\beta -1\right)  \left( x+\beta-
\alpha-3 \right)  \left( x+\alpha+\beta+2\,N+1 \right) }{ 4\left( x-1
 \right)  \left( x-3 \right) }}
$$
$$
G_1(x)={\frac { \left( \alpha+\beta \right)  \left( \alpha+\beta+ 2+2\,N \right) \left( x+1+\alpha-\beta
 \right)  }{ \left( 1-x^2 \right)
 \left( x+3 \right)  }}
$$
$$
G_2(x)={\frac { \left( \alpha-\beta \right)  \left( x-\alpha-\beta-1
 \right)  \left( x+\alpha+\beta+2\,N+1 \right) }{ \left( 1-x^2 \right)
 \left( x-3 \right) }}
$$

The form \re{HR_op} of the difference equation is preferable
because we here have the operator $H$ acting  directly on the
argument of the polynomials. The operator $H$ belongs to the class
of Dunkl shift operators: it contains both simple shifts $T^j$ and
the reflection operator $R$. Moreover, the operator $H$ preserves
the space of polynomials: it transforms any polynomials of degree
$n$ into a polynomial of the same degree $n$. Operators of this
kind were considered in the theory of Bannai-Ito polynomials
\cite{BI_TVZ}. However, in contrast to the Bannai-Ito situation,
the operator $H$ in the present case is of  {\it second order},
while the Bannai-Ito polynomials are eigenfunctions of a
Dunkl-shift operator of the first order \cite{BI_TVZ}.
Equivalently, this means that for generic values of the parameters
$\alpha, \beta$, the dual -1 Hahn polynomials are eigenvectors of
a 5-diagonal matrix (while the Bannai-Ito polynomials are
eigenvectors of a 3-diagonal matrix).

\bigskip\bigskip
{\Large\bf Acknowledgments}
\bigskip

The authors would like to gratefully acknowledge the hospitality
extended to LV and AZ by Kyoto University and to ST and LV by the
Donetsk Institute for Physics and Technology  in the course of
this investigation. The research of LV is supported in part by a
research grant from the Natural Sciences and Engineering Research
Council (NSERC) of Canada.

\newpage

\bb{99}

\bi{BI} E. Bannai and T. Ito, {\it Algebraic Combinatorics I: Association Schemes}. 1984. Benjamin \& Cummings, Mento Park, CA.

\bi{Chi_Bol} T.~Chihara, {\it On kernel polynomials and related
systems}, Boll. Unione Mat. Ital., 3 Ser., {\bf 19} (1964),
451--459.

\bi{Chi} T. Chihara, {\it An Introduction to Orthogonal
Polynomials}, Gordon and Breach, NY, 1978.

\bi{JSJ} E.~I.~Jafarov , N.~I.~Stoilova  and J.~Van der Jeugt,
{\it Finite oscillator models: the Hahn oscillator}, ArXiv:
1101.5310.

\bi{KLS} R.~Koekoek, P.~Lesky, R.~Swarttouw, {\it Hypergeometric
Orthogonal Polynomials and Their Q-analogues}, Springer-Verlag,
2010.

\bi{Leonard} D.Leonard, {\it Orthogonal Polynomials, Duality and
Association Schemes}, SIAM J. Math. Anal. {\bf 13} (1982)
656--663.

%\bi{Mag} A.P. Magnus, {\it Painleve\'e-type differential equations
%for the recurrence coefficients of semi-classical orthogonal
%polynomials}, J. Comp. Appl. Math. {\bf 57} (1995), 215-237.

%\bi{Mar} P.Maroni, {\it Variations around classical orthogonal polynomials.
%Connected problems}, J.Comp.Appl.Math. {\bf 48} (1993), 133-155.

%\bi{Mass} D.Masson, {\it Difference equations, continued fractions, Jacobi
%matrices and orthogonal polynomials}, in Nonlinear Numerical methods and
%Rational Approximations, (A.Cuyt, ed.), 239-257. Reidel Publ.Co. 1988.

\bi{NSU} A.F. Nikiforov, S.K. Suslov, and V.B. Uvarov, {\em
Classical Orthogonal Polynomials of a Discrete Variable},
Springer, Berlin, 1991.

%\bi{Pas} P.I. Pastro, {\it Orthogonal polynomials and some q-beta
%integrals of Ramanujan},\\  J.Math.Anal.Appl. {\bf 112} (1985),
%517--540.

\bi{Rosen} M. Rosenblum, {\it Generalized Hermite Polynomials and
the Bose-like Oscillator Calculus}, in: Oper. Theory Adv. Appl.,
vol. {\bf 73}, Birkhauser, Basel, 1994, pp. 369--396.
ArXiv:math/9307224.

%\bi{SVZ} V.Spiridonov, L.Vinet and A.Zhedanov, {\it Spectral
%transformations, self-similar reductions and orthogonal polynomials}, J.Phys.
%A:  Math.  and Gen.  {\bf 30} (1997), 7621--7637.

%\bi{Sz} G. Szeg\H{o}, Orthogonal Polynomials, fourth edition,  AMS, 1975.

%\bi{Ter} P. Terwilliger, {\it Two linear transformations each tridiagonal with respect to an eigenbasis of the other}.
%Linear Algebra Appl. {\bf 330} (2001) 149--203.

%\bi{VYZ} L.Vinet, O.Yermolayeva and A.Zhedanov, {\it A method to
%study the Krall and q-Krall polynomials}, J.Comp.Appl.Math. {\bf
%133} (2001) 647--656.

%\bi{VZ_Bochner} L.Vinet and A.Zhedanov, {\it Generalized Bochner
%theorem: characterization of the Askey-Wilson polynomials},
%J.Comp.Appl.Math., {\bf 211} (2008) 45 -- 56.

\bi{Stoilova} N.~I.~Stoilova  and J.~Van der Jeugt, {\it An
Exactly Solvable Spin Chain Related to Hahn Polynomials}, SIGMA
{\bf 7} (2011), 03.

\bi{Ter} P.~Terwilliger, {\it Two linear transformations each
tridiagonal with respect to an eigenbasis of the other}. Linear
Algebra Appl. {\bf 330} (2001) 149--203.

\bi{BI_TVZ} S.Tsujimoto, L.Vinet and A.Zhedanov, {\it Dunkl-shift
operators and Bannai-Ito polynomials}, arXiv:1106.3512 .

%\bi{Vidunas} R. Vid\=unas, {\it Askey-Wilson relations and Leonard pairs},  Discrete Mathematics {\bf 308} (2008) 479-–495. arXiv:math/0511509v2.

\bi{VZ_little} L.~Vinet and A.~Zhedanov, {\it A ``missing`` family
of classical orthogonal polynomials}, J. Phys. A: Math. Theor.
{\bf 44} (2011) 085201,    arXiv:1011.1669v2.

\bi{VZ_big} L.~Vinet and A.~Zhedanov, {\it A limit $q=-1$ for big
q-Jacobi polynomials}, Trans.Amer.Math.Soc., to appear,
arXiv:1011.1429v3.

\bi{VZ_Bochner} L.~Vinet and A.~Zhedanov, {\it A Bochner theorem
for Dunkl polynomials}, SIGMA, {\bf 7} (2011), 020;
arXiv:1011.1457v3.

%\bi{WW} E.T. Whittacker, G.N. Watson, {\em A Course of Modern
%Analysis}, Cambridge, 1927.

%\bi{Zhe} A. S. Zhedanov. {\it "Hidden symmetry" of Askey-Wilson polynomials}, Teoret. Mat. Fiz.
%{\bf 89} (1991) 190--204. (English transl.: Theoret. and Math. Phys. {\bf 89} (1991), 1146--1157).

%\bi{ZheL} A.Zhedanov, {The "classical" Laurent biorthogonal
%polynomials}, J. Comp. Appl.Math. {\bf 98} (1998), 121--147.

\eb

\end{document}